\documentclass[11pt]{article}

\usepackage{amsmath,amsfonts,amssymb}

\def\R{{\bf R}}

\def\K{{\bf K}}
\def\eps{{\varepsilon}}

\def\bE {{\mathbb E}}
\def\bP {{\mathbb P}}
\def\bR {{\mathbb R}}

\def\<{{\langle}}
\def\>{{\rangle}}

\def\qed{{\hfill $\Box$ \bigskip}}

\newtheorem{theorem}{Theorem}[section]
\newtheorem{lemma}[theorem]{Lemma}

\newtheorem{proposition}[theorem]{Proposition}

\newtheorem{definition}[theorem]{Definition}

\begin{document}

\title{\vspace*{-1.5cm}
Small time asymptotics for Brownian motion with singular drift}

\author{ Zhen-Qing Chen$^{1, 2}$ \quad Shizan Fang$^{3}$ \quad
Tusheng Zhang$^{4, 5}$}

\footnotetext[1]{\, Department of Mathematics, University of
Washington, Seattle, WA 98195, USA. Email: zqchen@uw.edu}
\footnotetext[2]{\, School of  Mathematics and Statistics, Beijing Institute of Technology, Beijing, China}
\footnotetext[3]{\,  Department
of Mathematics, University of Bourgogne, France.  Email: shizan.fang@u-bourgogne.fr}
\footnotetext[4]{\, School of Mathematics, University of Manchester, Oxford Road, Manchester
M13 9PL, England, U.K. Email: tusheng.zhang@manchester.ac.uk}
\footnotetext[5]{\, School of Mathematics, University of Science and Technology of China, Hefei, China.}
\maketitle

\begin{abstract}
We establish a small time large deviation principle and a Varadhan type asymptotics for Brownian motion with singular drift
on $\bR^d$ with $d\geq 3$
whose infinitesimal generator is $\frac12 \Delta + \mu \cdot \nabla$, where each $\mu_i$ of $\mu = (\mu_1, \cdots, \mu_d)$ is a measure in
some suitable Kato class.
\end{abstract}


\medskip


\section{Introduction}

For non-divergence form elliptic operator ${\cal L}=\frac12 \sum_{i, j=1}^d a_{ij}(x) \frac{\partial^2}{\partial x_i \partial x_j}$ on $\bR^d$ with bounded, symmetric and uniformly elliptic diffusion matrix $A(x)=(a_{ij}(x))$ that is uniformly H\"older continuous,
Varadhan showed the following small time asymptotics for the heat kernel $p(t, x, y)$ of
${\cal L}$ in \cite{V1}
\begin{equation}\label{e:1.1}
\lim_{t \downarrow 0} t \log p (t, x, y) =-d (x, y)^2 /2 \quad \hbox{for } x, y \in \bR^d,
\end{equation}
where $d(x, y)$ is the Riemannian metric induced by $A(x)^{-1}$.
Later, \eqref{e:1.1} is extended by Norris \cite{N} to divergence form elliptic operator
$$
{\cal L} = \frac1{2 \rho (x)} \sum_{i,j=1}^d
\frac{\partial}{\partial x_i} \left( \rho (x) a_{ij}(x) \frac{\partial}{\partial x_j}
\right),
$$
 where $A(x)=(a_{ij}(x))$ is symmetric locally bounded and locally uniformly elliptic
   and $\rho (x)$ is a measurable function that is locally bounded between two positive constants.
In fact, the results in \cite{N} are more general,  allowing $\bR^d$ being
replaced by a Lipschitz manifold.  In \cite{AH}, Ariyoshi and Hino studied the integral version
of \eqref{e:1.1} and showed that for heat semigroup $\{P_t;  t\geq 0\}$ associated with any
symmetric local Dirichlet form $({\cal E}, {\cal F})$ on $L^2(E; m)$ with $\sigma$-finite measure $m$,
\begin{equation}\label{e:1.2}
\lim_{t \downarrow 0} t \log P_t(A_1, A_2) =-d (A_1, A_2)^2 /2
\end{equation}
for any Borel subsets $A_i\subset E$ with $0< m(A_i) <\infty$ for $i=1, 2$.
Here $P_t (A_1, A_2):=\int_{A_2} P_t {\bf 1}_{A_1} (x) m(dx)$ and $d(A_1, A_2)$ is the intrinsic distance between
$A_1$ and $A_2$ induced by the local Dirichlet form $({\cal E}, {\cal F})$.
The above result \eqref{e:1.2} has recently been extended in \cite{HM}
to lower order perturbation of symmetric strongly local Dirichlet forms.
For earlier and other related work, see the references in \cite{Zh, N, HR}.
\vskip 0.3cm
In this paper, we study pointwise asymptotic property \eqref{e:1.1}
for Brownian motion with singular drifts   in $\bR^d$  with $d\geq 3$; that is,
\begin{equation}\label{e:1.3}
dX_t=dW_t+dA_t \quad \hbox{with }  X_0=x,
\end{equation}
where $A$ is a continuous additive functional of $X$ having ``Revuz measure $\mu$".
Informally $\{X_t, t\geq 0\}$ is a diffusion process in $\bR^d$ with generator $\frac{1}{2}\Delta+ \mu\cdot \nabla$, where $\mu=(\mu_1,\cdots ,\mu_d)$ is a vector-valued signed measure on $\bR^d$ belonging to the Kato class $\K_{d, 1}$ to be introduced below.  When $\mu_i (dx)=b_i(x)dx$ for some function $b_i$, $X$ is a solution to the stochastic differential equation
\begin{equation}\label{e:1.4}
dX_t=dW_t+b(X_t)dt \quad \hbox{with }  X_0=x
\end{equation}

\begin{definition}
A signed measure $\nu$ on $\bR^d$ with $d\geq 3$
is said to be in the Kato class $\K_{d,k}$ (for $k=1,2$) if
$$
\lim_{r \downarrow 0}\sup_{x\in \bR^d}\int_{|x-y|\leq r}\frac{|\nu|(dy)}{|x-y|^{d-k}}=0,
$$
where $|\nu|$ is the total variation of $\nu$.
A measurable function $f$ on $\bR^d$ is said to be in the Kato class $\K_{d,k}$ (for $k=1,2$)
if $|f(x)| dx \in \K_{d, k}$.
\end{definition}

\medskip

 Clearly any bounded measurable functions are in the Kato class  $\K_{d,k}$ for $k=1,2$.
By H\"older inequality, it is easy to see that $L^p(\bR^d; dx)\subset \K_{d, 1}$ for $p>d$
and  $L^p(\bR^d; dx)\subset \K_{d, 2}$ for $p>d/2$.
But a measure in $\K_{d, 1}$ can be quite singular.
It is shown in  \cite[Proposition 2.1]{BC} that for a Borel measure $\mu$ on $\bR^d$,
if there are constants $\kappa >0$ and $\gamma >0$ so that
  $\mu (B(x, r) \leq \kappa r^{d-1+\gamma}$ for all $x\in \bR^d$ and $r\in (0, 1]$,
  then $\mu \in \K_{d, 1}$.  Thus in particular,  if $A\subset \bR^d$ is an Alfors $\lambda$-regular set
  with $\lambda \in (d-1, d]$, then the Hausdorff measure ${\mathcal  H}^\lambda$ restricted to
  $A$ is in $\K_{d, 1}$.

\medskip

To recall the precise definition of a Brownian motion with a measure drift  $\mu=(\mu_1,\cdots ,\mu_d)$ with $\mu_i \in \K_{d, 1}$
for $1\leq i\leq d$,
fix a non-negative smooth function $\varphi$ in $\bR^d$ with compact support and $\int\varphi(x)dx=1$.  For any positive integer $n$, we put $\varphi_n(x)=2^{nd}\varphi(2^nx)$. For $1\leq i\leq d$, define
$$
b^{(n)}_i(x)=\int\varphi_n(x-y)\mu_i(dy)$$
Put $b^{(n)}=(b^{(n)}_1, \cdots, b^{(n)}_d)$. The following definition is taken from \cite{BC}.

\begin{definition}\rm
A Brownian motion with drift $\mu$ is a family of probability measures $\{\bP_x, x\in \bR^d\}$ on $C([0,\infty ), \bR^d)$, the space of continuous functions on $[0, \infty )$, such that under each $\bP_x$ the coordinator process $X$ has the decomposition
$$X_t=x+W_t+A_t$$
where
\begin{description}
\item{(a)} $A_t= (A^{(1)}_t\, \cdots, A^{(d)}_t) =
\lim_{n\rightarrow \infty}\int_0^t b^{(n)}(X_s)ds$ uniformly in $t$ over finite intervals, where the convergence is in probability;

\item{(b)} there exists a subsequence $\{n_k\}$ such that for every $t>0$,
$$
\sup_{k\geq 1} \int_0^t |b^{n_k}(X_s)| ds<\infty \quad a.s.
$$

\item{(c)} $W_t$ is a standard Brownian motion in $\bR^d$ starting from the origin.
\end{description}
\end{definition}

It is established in \cite{BC} that, when each $\mu_i$ is in the Kato class $\K_{d, 1}$, Brownian motion
with drift  $\mu=(\mu_1, \cdots, \mu_d)$, denoted by $\{X_t, t\geq 0\}$,  exists and is unique in law for every starting point $x\in \bR^d$.
Moreover, it is shown there  that
 $X=\{X_t, t\geq 0; \bP_x, x\in \bR^d\}$ forms a conservative Feller process on $\bR^d$.
Later,  it is obtained in \cite{KS}  that   $X$ has a jointly continuous
transition density function $q(t,x,y)$ which admits the following Gaussian type estimate:
 \begin{equation}\label{1.2-1}
 C_1e^{-C_2 t}t^{-{d}/{2}}e^{- {C_3|x-y|^2}/{t}}\leq q(t,x,y)\leq C_4e^{C_5t}t^{- {d}/{2}}e^{-{C_6|x-y|^2}/{t}}
 \end{equation}
for all $(t,x,y)\in (0, \infty)\times \bR^d\times \bR^d$, where $C_i, 1\leq i\leq 6$,  are some positive constants.

\medskip

In this paper, we are concerned with the  precise Varadhan type small time asymptotics of the transition density
function $q(t, x, y)$ of $X$:
 \begin{equation}\label{1.3}
\lim_{t\rightarrow 0}t\log q(t,x,y)=-\frac{|x-y|^2}{2}
\end{equation}
for every $x, y\in \bR^d$.
 Note that since the drift measure $\mu$ is merely assumed to be in  Kato class  $\K_{d, 1}$, the law of $X$ can be singular
 with respect to the law of Brownian motion. Moreover, it is not clear if Browian motion with drift $\mu$ can be constructed
 via sectorial non-symmetric Dirichlet form (see \cite[p.793]{BC}).
 One can not even apply the result from \cite{HM} to
 get the integral version of small time asymptotics result \eqref{e:1.2}
 for our Brownian motion with drift $\mu$.  Hence a new approach is needed.
 %end new

For each $x\in \bR^d$, denote by  $\nu^x_{\varepsilon}$   the law of $\{X_{\varepsilon t}, 0\leq t\leq 1\}$ on $C([0,1],\bR^d)$ under $\bP_x$.
It is known (see \cite{V1, V2}) that the Varadhan small time asymptotics (\ref{1.3}) is closely related to the small time large deviation principle of the family $\{\nu^x_{\varepsilon}, \varepsilon >0\}$.
In fact, in this paper,    we will establish  a small time large deviation principle for $\{ \nu^x_{\varepsilon}; \eps>0\}$ for every $x\in \bR^d$
 by showing that $\{ \nu^x_{\varepsilon}; \eps>0\}$ is exponentially equivalent to the corresponding laws of Brownian motion.
 This  is done in Section \ref{S:2}. The upper bound  small time asymptotics of \eqref{1.3} is obtained by the construction of heat kernel of $X$ and their pointwise upper bound estimates. This is carried out in Proposition \ref{P:3.3}.
 To establish the lower bound small time asymptotics of \eqref{1.3} of transition density $q(t,x,y)$, we will use a crucial estimate ( Proposition 2.2) of the Laplace transform of the  drift $A$ to get a lower bound of the probability that the process $X$ belongs to a small ball. This is the content of Proposition \ref{P:3.4} of Section \ref{S:3}.

\section{Small time large deviations}\label{S:2}
\setcounter{equation}{0}

In this section, we will establish a small time large deviation principle for   Brownian motion with measure drift $\mu$.
For a signed measure $\mu$ on  $ \bR^d$, define
\begin{eqnarray}\label{e:2.1}
N_t^{\alpha}(\mu)&:=&\sup_{x \in \bR^d }\int_0^t\int_{\bR^d}s^{- {(d+1)}/{2}}
\exp \left(-\alpha\frac{|x-y|^2}{s} \right)|\mu|(dy)ds,
\end{eqnarray}
where $|\mu|$ denotes the measure of total variation.
Recall (see e.g. \cite{KS}) that $\mu\in \K_{d,1}$ if and only if $\lim_{t\rightarrow 0}N_t^{\alpha}(\mu)=0$ for any $\alpha>0$. For the transition kernel $q(t,x,y)$ and a
 signed measure $\mu$ on  $ \bR^d$, we also introduce
 \begin{eqnarray}\label{2.0-1}
\Lambda_t(\mu)&:=&\sup_{x \in \bR^d}\int_0^t\int_{\bR^d}\frac{q(s,x,y)}{\sqrt{s}}|\mu|(dy)ds.
\end{eqnarray}
As a consequence of the two-sided Gaussian bounds of the transition density (\ref{1.2-1}), it is easy to see that $\mu$ is in the Kato class $\K_{d,1}$ if and only if $\lim_{t\rightarrow 0}\Lambda_t(\mu)=0$. For a function $f$ in $\K_{d, 1}$, we write $\Lambda_t (f)$ for $\Lambda_t (fdx)$.
 \vskip 0.3cm
 Let $X$ be the Brownian motion with measure drift $\mu$ stated in Section 1 and  $\bP_x$ the distribution of $X$ starting from $x$.
 For $x\in \bR^d$ and $\eps>0$, denote by  $\nu^x_{\varepsilon}$ the distribution of $\{X_{\eps t}; t\in [0, 1]\}$ under $\bP_x$ on the path space $C([0,1],\bR^d)$.    We have the following large deviation result.

\begin{theorem}\label{T:2.1} For each fixed $x\in \bR^d$,
$\{\nu^x_{\varepsilon}, \varepsilon>0\}$ satisfies a large deviation principle with rate function
$$ I(f)=\left\{\begin{array}{ll} \frac{1}{2}\int_0^1|\dot{f}(t)|^2dt & \mbox{if $f$ is absolutely continuous},\\
 \infty& \mbox{otherwise, }
 \end{array}\right.
 $$
 where $\dot{f}$ stands for the derivative of $t\rightarrow f(t)$.
 Namely,
 \begin{description}
\item{\rm (i)} for any closed subset $C\subset C ([0,1], \bR^d)$,
 $$\limsup_{\varepsilon\rightarrow 0}\varepsilon \log \nu^x_{\varepsilon}(C)\leq -\inf_{f\in C (x)}I(f),$$
where $C(x ) =\{ f\in C: f(0)=x\}$;

\item{\rm  (ii)} for any open subset $G\subset C([0,1], \bR^d)$,
 $$\liminf_{\varepsilon\rightarrow 0}\varepsilon \log \nu^x_{\varepsilon}(G)\geq -\inf_{f\in G (x)}I(f),
 $$
where $G(x)=\{f\in G: f(0)=x\}$.
\end{description}
\end{theorem}

\medskip

To prove  Theorem \ref{T:2.1}, we need the following crucial estimate for the Laplace transform of an additive functional.

\begin{proposition}\label{P:2.2}
Let $b(\cdot)$ be a positive function in the Kato class $\K_{d, 1}$. Then, for any $\lambda>0$
\begin{equation}\label{2.0-2}
 \bE_x \left[ e^{\lambda \int_0^tb(X_s)ds} \right]
 \leq  \left(1+\lambda \sqrt{t}\Lambda_t(b) \right) \exp \left( \lambda^2 t\Lambda_t(b)^2 \right)
 %new
 \leq 2 \exp \left( 2\lambda^2 t\Lambda_t(b)^2 \right),
\end{equation}
where $\Lambda_t(b)$ is defined as in  \eqref{2.0-1}.
\end{proposition}

\noindent{\bf Proof.} We claim that for non-negative integer $n\geq 0$,
\begin{equation}\label{2.0-3}
 \bE_x \left[ \Big(\int_0^tb(X_s)ds \Big)^n \right]\leq n! \, \alpha_n  \left( \sqrt{t} \Lambda_t(b) \right)^n ,
\end{equation}
where $\alpha_0=\alpha_1=1$, $\alpha_n=\Pi_{k=2}^n  \left(1-\frac{1}{k}\right)^{(k-1)/{2}}
 \left( \frac{1}{k} \right)^{1/2} $ for $n\geq 2$.
Indeed, for $n\geq 1$, we have by the Markov property of $X$,
\begin{eqnarray}\label{2.0-4}
 &&\bE_x \left[ \Big( \int_0^tb(X_s)ds \Big)^n \right] \nonumber\\
 &=&n!\int_0^t ds_n\int_0^{s_n}ds_{n-1}\cdots \int_0^{s_{2}}ds_1\bE_x \left[b(X_{s_1})b(X_{s_{2}})\cdots b(X_{s_n})
 \right]  \nonumber \\
 &=&n!\int_{\{0\leq  s_1<s_2< \cdots <s_n \leq t\}}  ds_1 \cdots ds_n  \int_{(\bR^d)^{\otimes n}} b(y_1)\cdots b(y_n)q(s_1,x,y_1)q(s_2-s_1, y_1, y_2)\nonumber\\
 &&\quad\quad\quad\cdots q(s_n-s_{n-1}, y_{n-1}, y_n)dy_n\cdots dy_1.
\end{eqnarray}
When $n=1$,
\begin{eqnarray}\label{e:2.6}
\sup_{x\in \bR^d} \bE_x \left[ \int_0^tb(X_s)ds   \right]
&=& \sup_{x\in \bR^d} \int_0^t \int_{\bR^d}  b(y) q(s, x, y)  dy  ds  \nonumber\\
&\leq & \sqrt{t} \sup_{x\in \bR^d} \int_0^t \int_{\bR^d}  b(y) \frac{q(s, x, y)}{\sqrt{s}}  dy  \nonumber\\
&=& \sqrt{t} \Lambda_t (b).
\end{eqnarray}
So \eqref{2.0-3} holds for $n=0$ and $n=1$ as claimed.
When $n=2$, by \eqref{2.0-4} and \eqref{e:2.6},
\begin{eqnarray*}
&& \sup_{x\in \bR^d} \bE_x \left[ \Big( \int_0^tb(X_s)ds \Big)^2  \right] \\
&=& 2 \sup_{x\in \bR^d} \int_0^t \int_{\bR^d}  b(y_1) q(s_1, x, y_1) \left( \int_{s_1}^t \int_{\bR^d}  b(y_2) q(s_2-s_1, y_1, y_2) ds_2 dy_2) \right)  ds_1 dy_1 \\
&=&  2 \sup_{x\in \bR^d} \int_0^t \int_{\R^d}  b(y_1) q(s_1, x, y_1) \left( \int_0^{t-s_1} \int_{\bR^d}  b(y_2) q(r, y_1, y_2) dr dy_2) \right)  ds_1 dy_1 \\
 &\leq &    2 \sup_{x\in \bR^d} \int_0^t \int_{\bR^d}  b(y_1) q(s_1, x, y_1)  \sqrt{t-s_1} \Lambda_t (b)ds_1 dy_1 \\
&= & 2 \Lambda_t (b)  \sup_{x\in \bR^d} \int_0^t \int_{\bR^d}  \sqrt{(t-s_1)s_1} \, b(y_1) \frac{q(s_1, x, y_1)}{  \sqrt{s_1}}  ds_1 dy_1 \\
&\leq& 2! \cdot \frac12\cdot \left(\sqrt{t} \Lambda_t (b)\right)^2,
\end{eqnarray*}
as $\max_{s_1 \in [0, t] } \sqrt{ s_1(t-s_1)} = t/2$. This shows that \eqref{2.0-3} holds for $n=2$.
Now assuming \eqref{2.0-3} holds for $n=k\geq 2$, we next show it holds for $n=k+1$.
By \eqref{2.0-3} for $n=k$,
\begin{eqnarray*}
&& \sup_{x\in \bR^d}  \bE_x \left[ \Big( \int_0^tb(X_s)ds \Big)^{k+1} \right] \nonumber\\
&\leq& (k+1) \cdot k! \int_0^t \sup_{x\in \bR^d }
\int_{\bR^d}  b(y_1)  q(s_1,x,y_1)\left( \int_{\{0<s_2-s_1< \cdots <s_k  -s_1\leq t-s_1\}}  ds_s \cdots ds_n   \right. \nonumber\\
 && \left. \quad   \int_{(\bR^d)^{\otimes k}}  b(y_2)\cdots b(y_k)q(s_2-s_1, y_1, y_2) \cdots q(s_n-s_{n-1}, y_{n-1}, y_k) dy_k \right) ds_1 dy_1\\
 &\leq & (k+1) ! \alpha_k   \Lambda_t (b)^k  \sup_{x\in \bR^d}  \int_0^t \int_{\bR^d} b(y_1)  \frac{ q(s_1, x, y_1)}{\sqrt{s_1}} s_1^{1/2} (t-s_1)^{k/2} ds_1 dy_1  \\
 &\leq &  (k+1) ! \alpha_k \left(\frac{k}{k+1}\right)^{k/2}\left(\frac{1}{k+1}\right)^{1/2} t^{(k+1)/2} \Lambda_t (b)^{k+1} \\
 &=& (k+1) ! \alpha_{k+1} \left( \sqrt{t} \Lambda_t (b)\right)^{k+1},
\end{eqnarray*}
as $\max_{s_1\in [0, s_1]}  s_1^{1/2} (t-s_1)^{k/2} = \left(\frac{k}{k+1}\right)^{k/2}  \left(\frac{1}{k+1}\right)^{1/2} t^{(k+1)/2}$.
This shows that \eqref{2.0-3} holds for $n=k+1$. By induction,  we have established that
 \eqref{2.0-3} holds for all $n\geq 1$.

 \smallskip

Observe that $\alpha_n\leq \frac{1}{\sqrt{n!}}$. We have by  \eqref{2.0-3}  that for $\lambda>0$,
\begin{equation}\label{2.0-17}
 \bE_x \left[e^{\lambda \int_0^tb(X_s)ds}\right]
 = \sum_{n=0}^\infty \frac1{n!} \bE_x \left[  \Big( \lambda \int_0^tb(X_s)ds\Big)^n \right]
 \leq \sum_{n=0}^{\infty}\frac{ \left(\lambda \sqrt{t} \Lambda_t(b) \right)^n }{\sqrt{n!}}.
\end{equation}
Since  $(2n+1)! \geq (2n)!\geq (n!)^2$  for integer $n\geq 0$, we have for $z\geq 0$,
 \begin{eqnarray}\label{2.0-18}
 \Phi(z) &:=&\sum_{n=0}^{\infty}\frac{z^n}{\sqrt{n!}}
  =  \sum_{n=0}^{\infty}\frac{z^{(2n)}}{\sqrt{(2n)!}}
+\sum_{n=0}^{\infty}\frac{z^{2n+1}}{\sqrt{(2n+1)!}} \nonumber \\
 &\leq &\sum_{n=0}^{\infty}\frac{z^{2n}}{n!}
+z\sum_{n=0}^{\infty}\frac{z^{2n}}{n!}
 =  (1+z)e^{z^2}.
\end{eqnarray}
This combined with \eqref{2.0-17} yields that
\begin{equation}\label{2.0-19}
 \bE_x \left[ e^{\lambda \int_0^tb(X_s)ds} \right]
 \leq  \left(1+\lambda\Lambda_t(b)\sqrt{t} \right) e^{\lambda^2\Lambda_t(b)^2t}.
\end{equation}
This completes the proof of the proposition as
$ 1+a  \leq 2 e^{ a^2}$ for every $a\geq 0$.
\qed

\vskip 0.4cm

\noindent{\bf Proof of Theorem \ref{T:2.1}.}   Fix $x\in \bR^d$. Let $\mu_{\varepsilon}^x$ denote the law of the Brownian motion $\{x+W_{\eps t}; t\in [0, 1]\}$ on the path space $C([0,1],\bR^d)$. Then it is well known that $\{\mu_{\varepsilon}^x, \varepsilon>0\}$ obeys a large deviation principle with a rate function $I(\cdot)$ defined as in the statement of Theorem 2.1. Note that
 $$X_{\eps t}=x+W_{\eps t}+A_{\eps t}.$$
 According to \cite[Theorem 4.2.13]{DZ}, to prove that $\nu_{\varepsilon}^x$ satisfies a large deviation principle with the same rate function as $\mu_{\varepsilon}^x$ it is sufficient to show that the two families $\{\nu_{\varepsilon}^x\}$, $\{\mu_{\varepsilon}^x\}$ are exponentially equivalent, namely for any $\delta>0$,
\begin{equation}\label{2.0-20}
 \limsup_{\varepsilon\rightarrow 0}\varepsilon \log \bP_x
 \left( \sup_{0\leq t\leq 1}|A_{\eps t}|=
 \sup_{0\leq t\leq 1}|X_{\eps t}-(x+W_{\eps t})|>\delta \right)=-\infty.
\end{equation}
To this end, we first deduce an estimate similar to (\ref{2.0-2}) for the Laplace transform of the process $A$, namely, for any $\lambda>0$
\begin{eqnarray}\label{2.0-21}
   \bE_x \left[e^{\lambda \sup_{0\leq s\leq t}|A_s|} \right]
    \leq2 \exp\left( 2\lambda^2 t  C_4^2 e^{2 C_5}N^{C_6}_t \left( \hbox{$\sum_{i=1}^d$} |\mu_i|\right)^2 \right) .
 \end{eqnarray}
For $1\leq i\leq d$, recall the following function defined in Section 1
$$
b^{(n)}_i(x)=\int\varphi_n(x-y)\mu_i(dy).$$
We claim that
\begin{equation}\label{2.1}
 \Lambda_t(|b^{(n)}_i|)\leq C_4e^{C_5}N^{C_6}_t(|\mu_i|),
\end{equation}
where $N^{C_6}_t(|\mu_i|)$ was defined in \eqref{e:2.1}.
Note that
$$
|b^{(n)}_i|(x)\leq \int\varphi_n(x-y)|\mu_i|(dy).$$
Using the upper  Gaussian-type estimate  \eqref{1.2-1} for $q(s,x,y)$,  we have for $t\leq 1$,
\begin{eqnarray*}
 && \Lambda_t(|b^{(n)}_i|) \\
 &\leq &\sup_{x}\int_0^t\int_{\bR^d}\frac{q(s,x,y)}{\sqrt{s}}
 \left(\int_{\bR^d}\varphi_n(y-z)|\mu_i|(dz) \right)dy ds\nonumber\\
 &\leq& \sup_{x}\int_0^t\int_{\bR^d}|\mu_i|(dz)\frac{1}{\sqrt{s}}
 \left( \int_{\bR^d}q(s,x,y^{\prime}+z)\varphi_n(y^{\prime})|dy^{\prime} \right) ds\nonumber\\
 &\leq& \sup_{x}\int_0^t\int_{\bR^d}|\mu_i|(dz)\frac{1}{\sqrt{s}}
 \left(\int_{\bR^d} C_4e^{C_5s}s^{- {d}/{2}}e^{-{C_6|x-y^{\prime}-z|^2}/{s}}\varphi_n(y^{\prime})|dy^{\prime} \right) ds\nonumber\\
 &\leq&C_4e^{C_5}\int_{\bR^d}\varphi_n(y^{\prime})|dy^{\prime}\sup_{x, y^{\prime}}\int_0^t
 \left(\int_{\bR^d}s^{- {d+1}/{2}}e^{-{C_6|x-y^{\prime}-z|^2}/{s}}|\mu_i|(dz) \right) ds\nonumber\\
 &\leq& C_4e^{C_5}N^{C_6}_t(|\mu_i|),
\end{eqnarray*}
where the fact $\int_{\bR^d}\varphi_n(y^{\prime})|dy^{\prime}=1$ was used in the last inequality.
This proves the claim \eqref{2.1}.
 We have by the definition of the process $A$,   Fatou's Lemma, Proposition 2.2 and  \eqref{2.1} that
 for every $\lambda >0$,
\begin{eqnarray*}
 && \bE_x \left[e^{\lambda \sup_{0\leq s\leq t}|A_s|} \right]   \\
 &\leq &\bE_x[e^{\lambda \sup_{0\leq s\leq t}\sum_{i=1}^d |A^{(i)}_s|}]\nonumber\\
   &\leq&\liminf_{n\rightarrow \infty}\bE_x[e^{\lambda  \int_0^t  \sum_{i=1}^d |b^{(n)}_i|(X_s) ds}]\nonumber\\
    &\leq& 2
    \exp \left( 2 \lambda^2 t\Lambda_t \left( \hbox{$\sum_{i=1}^d$} |b^{(n)}_i| \right)^2 \right)\nonumber\\
 &\leq& 2
   \exp\left( 2\lambda^2 t  C_4^2 e^{2 C_5}N^{C_6}_t \left( \hbox{$\sum_{i=1}^d$} |\mu_i|\right)^2 \right) .
   \end{eqnarray*}
 This establishes the claim  (\ref{2.0-21}). We are now ready to prove (\ref{2.0-20}).
 For any $\delta, \lambda>0$,  we have by  (\ref{2.0-21}),
\begin{eqnarray}\label{2.5}
 &&\bP_x \left(\sup_{0\leq t\leq 1}|A_{\eps t}|>\delta \right)
 =\bP_x \left( \sup_{0\leq s\leq \varepsilon}|A_{s}|>\delta \right)\nonumber\\
 &\leq& e^{-\lambda \delta}  \bE_x \left[ e^{\lambda \sup_{0\leq s\leq \varepsilon}|A_s|} \right]\nonumber\\
  &\leq&  2 e^{-\lambda \delta}
   \exp\left( 2\lambda^2 \eps  C_4^2 e^{2 C_5}N^{C_6}_\eps \left( \hbox{$\sum_{i=1}^d$} |\mu_i|\right)^2 \right) .
  \end{eqnarray}
Taking
$$
 \lambda =\frac{\delta}{4 \eps  C_4^2 e^{2 C_5}N^{C_6}_\eps \left( \hbox{$\sum_{i=1}^d$} |\mu_i|\right)^2}
$$
yields
\begin{eqnarray}\label{2.6}
    \bP_x \left(\sup_{0\leq t\leq 1}|A_{\eps t}|>\delta\right)
   \leq 2  \exp \left(  -\frac{\delta^2}{8\eps  C_4^2 e^{2 C_5}N^{C_6}_\eps \left( \hbox{$\sum_{i=1}^d$} |\mu_i|\right)^2 }\right)  .
\end{eqnarray}
Hence
\begin{eqnarray*}\label{2.7}
 && \limsup_{\eps \to 0} \varepsilon\log \bP_x \left( \sup_{0\leq t\leq 1}|A_{\eps t}|>\delta \right)  \\
 &\leq&  \limsup_{\eps \to 0} \left( \varepsilon\log  2 -  \frac{\delta^2}{8  C_4^2 e^{2 C_5}N^{C_6}_\eps
 \left( \hbox{$\sum_{i=1}^d$} |\mu_i|\right)^2 } \right)\\
 &=& -\infty,
 \end{eqnarray*}
where the last inequality is due to the fact that  $ \lim_{\eps \to 0}N^{C_6}_\eps \left( \hbox{$\sum_{i=1}^d$} |\mu_i|\right) = 0$.
  This completes the  proof of Theorem \ref{T:2.1} . \qed

\section{Small time asymptotics}\label{S:3}
\setcounter{equation}{0}

Recall that $q(t,x,y)$ denotes the transition density of the Markov process $\{X_t,t\geq 0\}$.
In this section we aim to establish the following Varadhan's asymptotics.

\begin{theorem}\label{T:3.1}
\begin{equation}\label{3.1}
\lim_{t\rightarrow 0}t\log q(t,x,y)=-\frac{|x-y|^2}{2}
\end{equation}
uniformly in   $x, y\in \bR^d$ such that $|x-y|$ is bounded.
 \end{theorem}

\medskip

We will establish the upper bound and the lower bound separately.
 For $a>0$, let
$$
G_a(s,x,y)=s^{-{d}/{2}}\exp  \left(- \frac{a |x-y|^2}{2s} \right).
$$
 Recall the following result from \cite[Lemma 3.1]{Z1}.

\begin{lemma}\label{L:3.2}
Suppose $0<a_1<a_2$. There exist positive constants $C_0 $ and $\alpha $ only depending on $a_1,a_2$ such that
$$\int_0^t\int_{\bR^d}G_{a_1}(t-s,x,z)|b(z)|\frac{G_{a_2}(s,z,y)}{s^{1/2}}dzds
\leq C_0 N_t^{\alpha }(|b|)G_{a_1}(t,x,y).$$
\end{lemma}

\medskip

\begin{proposition}\label{P:3.3}
Assume the drift measure $\mu =(\mu_1, \cdots, \mu_d)$ belongs to $\K_{d,1}$. Then
 $$
\limsup_{t\rightarrow 0}t\log q(t,x,y)\leq -\frac{|x-y|^2}{2}
$$
 uniformly in  $x, y\in \bR^d$ such that $|x-y|$ is bounded.
\end{proposition}

\noindent{\bf Proof}.
We will prove the proposition by showing that for any $\delta>0$, there exist constants $T_{\delta}>0$, $C_{2,\delta}$ such that for $t\in (0, T_{\delta}],$
\begin{equation}\label{3.2-2}
 q(t,x,y)\leq C_{2,\delta}t^{- {d}/{2}}\exp \left(-(1-\delta)\frac{|x-y|^2}{2t} \right) .
\end{equation}
Using the approximation arguments as in \cite{KS}, it suffices to show that (\ref{3.2-2}) holds for drift measure $\mu$ given by $\mu(dy)=b(y)dy$ and that the constant $T_{\delta}$ is determined only by the quantity $N_t^{\alpha}(\mu)$. From now on we suppose
$\mu(dy)=b(y)dy$ for some function $b=(b_1,...,b_d)\in \K_{d, 1}$.
Note that for $0<\delta <1$, we have
\begin{eqnarray}
|\nabla_z G_1(s,z,y)|
&=& \frac{|z-y|}{s}  \exp \left( -\frac{\delta}{2} \frac{|z-y|^2}{2s} \right) G_{(1-\frac{\delta}{2})} (s, z, y)  \nonumber \\
&\leq&  m_{\delta} s^{-1/2} G_{(1-\frac{\delta}{2})} (s, z, y), \label{e:3.26}
\end{eqnarray}
where  $m_{\delta }:=\sup_{r>0} r e^{- \delta r^2/2}= 1/\sqrt{e\delta}$.
Hence  for every $0<\delta<1$, by Lemma \ref{L:3.2} (with  $a_1=1-\delta$ and  $a_2=1- (\delta/2)$),
there are positive constants $C_\delta$ and $c_\delta$
depending only on $\delta$ so that
\begin{eqnarray}\label{3.3}
 J(t,x,y)&:=& \int_0^t\int_{\bR^d}G_{(1-\delta)}(t-s,x,z)|b(z)| |\nabla_z G_1(s,z,y)| \, dzds\nonumber\\
&\leq& C_{\delta} N_t^{c_{\delta}}(|b|)G_{(1-\delta)}(t,x,y).
\end{eqnarray}

Define recursively $I_k(t,x,y)$ as follows:
\begin{eqnarray*}
I_0(t,x,y)&=& p(t,x,y)=(2\pi t)^{-{d}/{2}}\exp \left(-\frac{|x-y|^2}{2t} \right), \\
 I_{k+1}(t,x,y) &=&\int_0^t\int_{\bR^d}I_{k}(t-s,x,z)b(z)\cdot \nabla_z p(s,z,y)dzds
 \quad \hbox{for } k\geq 0 .
\end{eqnarray*}
Recall the following identity (see, e.g., \cite{KS, Z1}):
\begin{equation}\label{3.9}
q(t,x,y)=\sum_{k=0}^{\infty}I_k(t,x,y).
\end{equation}
We next  show that for any $\delta \in (0, 1)$, it holds that
 \begin{equation}\label{3.10}
|I_k(t,x,y)|\leq (C_{\delta}N_t^{c_{\delta}}(|b|))^k t^{-{d}/{2}}\exp
\left(-(1-\delta)\frac{|x-y|^2}{2t} \right).
\end{equation}
In view of (\ref{3.3}), (\ref{3.10}) clearly holds for $k=0,1$. Suppose (\ref{3.10}) holds for $k\geq 1$.
By induction and \eqref{3.3}, we  have
\begin{eqnarray*}\label{3.11}
I_{k+1}(t,x,y)&\leq &(C_{\delta}N_t^{ {\delta}/{2}}(|b|))^k\int_{0}^{t} \int_{\bR^d}
G_{1-\delta}(t-s, x, z) |b(z)| |\nabla_z G_1 (s, z, y) | dz ds \nonumber\\
&\leq & (C_{\delta}N_t^{c{\delta}}(|b|))^k C_{\delta}N_t^{c_{\delta}}(|b|) G_{1-\delta} (t, x, y)\nonumber\\
&=&(C_{\delta}N_t^{c{\delta}}(|b|))^{k+1}t^{-{d}/{2}}\exp \left(-(1-\delta)\frac{|x-y|^2}{2t}\right).
\end{eqnarray*}
Since $\mu= b(x) dx\in \K_{d,1}$, there exists a constant $T_{\delta}>0$ such that
$C_{\delta}N_t^{{\delta}/{2}}(|b|)\leq \frac{1}{2}$ for   $t\leq T_{\delta}$. This together with (\ref{3.9}), (\ref{3.10}) gives that
\begin{equation*}\label{3.12}
q(t,x,y)\leq 2 t^{-{d}/{2}}\exp \left(-(1-\delta)\frac{|x-y|^2}{2t}\right).
\end{equation*}
From the proof above, we see that the constant $T_{\delta}$ only depends on the rate at which
 $N_t^{{\delta}/{2}}(|b|)$ tends to zero.  Consequently,
 $$
 t\log q (t, x, y) \leq - (1-\delta) \frac{|x-y|^2}{2} - (d/2)  t \log (2t) \qquad \hbox{for } t\in (0, T_\delta].
 $$
 It follows that $\limsup_{t\to 0}t\log q (t, x, y) \leq -\frac{|x-y|^2}{2}$ uniformly on compact subsets of $\bR^d\times \bR^d$,
 By an approximation procedure as in \cite{KS}, we asset that the Proposition hold also for $\mu \in \K_{d,1}$.
 \qed

\begin{proposition}\label{P:3.4}
Suppose $\mu=(\mu_1, \cdots, \mu_d) \in \K_{d,1}$. Then
\begin{equation}\label{3.13}
\liminf_{t\rightarrow 0}t\log q(t,x,y)\geq -\frac{|x-y|^2}{2}
\end{equation}
  uniformly in  $x, y\in \bR^d$ such that $|x-y|$ is bounded.
\end{proposition}

\noindent {\bf Proof}.
 %We will use a trick from  the proof of Theorem 4.12. in \cite{V1}.
 For $\varepsilon>0$, set $B(y,\varepsilon)=\{z; |z-y|<\varepsilon\}$. Let $r>0$. We first like to give an estimate
 for the probability $\bP_x(X_r\in B(y,\varepsilon))$. Let $\delta\in (0, \eps)$.
 We have
 \begin{eqnarray}\label{3.14}
&&\bP_x(W_r+x\in B(y,\varepsilon-\delta)) =\bP_x(X_r-A_r\in B(y,\varepsilon-\delta))\nonumber\\
&\leq & \bP_x(X_r-A_r\in B(y,\varepsilon-\delta), |A_r|<\delta)+\bP_x( |A_r|\geq \delta)\nonumber\\
&\leq& \bP_x(X_r\in B(y,\varepsilon))+\bP_x( |A_r|\geq \delta).
\end{eqnarray}
In view of (\ref{2.6}),
\begin{eqnarray}\label{3.15}
\bP_x( |A_r|\geq \delta) \leq 2  \exp \left(  -\frac{\delta^2}{8r  C_4^2 e^{2 C_5}N^{C_6}_r \big( \hbox{$\sum_{i=1}^d$} |\mu_i|\big)^2 }\right).
\end{eqnarray}
On the other hand,
\begin{eqnarray}\label{3.16}
&&\bP_x(W_r+x\in B(y,\varepsilon-\delta))=\int_{B(0,\varepsilon-\delta)}(2\pi r)^{-\frac{d}{2}}
e^{-\frac{|z-(y-x)|^2}{2r}}dz\nonumber\\
&\geq& (2\pi r)^{-{d}/{2}} \omega_d (\eps-\delta)^d \exp\left(-\frac{(|x-y|+\varepsilon-\delta)^2}{2r}\right),
\end{eqnarray}
where $\omega_d$ is the volume of the unit ball in $\R^d$.
Thus for $0<\delta<\varepsilon$,
\begin{eqnarray}\label{3.17}
\bP_x(X_r\in B(y,\varepsilon))
&\geq& (2\pi r)^{-d/2} \omega_d (\eps-\delta)^d \exp \left( -\frac{(|x-y|+\varepsilon-\delta)^2}{2r}\right)
\nonumber\\
&& -2  \exp \left(  -\frac{\delta^2}{8r  C_4^2 e^{2 C_5}N^{C_6}_r \big( \hbox{$\sum_{i=1}^d$} |\mu_i|\big)^2 }\right).
\end{eqnarray}
 Note for every $0<\eta<1$ and $\varepsilon>0$,
\begin{eqnarray}\label{3.18}
q(t,x,y)&= &\int_{\bR^d}q((1-\eta)t,x,z)q(\eta t,z,y)dz \nonumber\\
&\geq & \int_{B(y,\varepsilon)}q((1-\eta)t,x,z)q(\eta t,z,y)dz\nonumber\\
&\geq& \inf_{z\in B(y,\varepsilon)}q(\eta t,z,y) \int_{B(y,\varepsilon)}q((1-\eta)t,x,z)dz \nonumber\\
&=&\inf_{z\in B(y,\varepsilon)}q(\eta t,z,y) \bP_x \left( X_{(1-\eta)t}\in B(y,\varepsilon) \right).
\end{eqnarray}

Let $M>0$. We have by  \eqref{3.18}, \eqref{1.2-1}  and \eqref{3.17}) with $r=(1-\eta)t$
that for any $x, y\in \R^d$ with $|x-y|\leq M$,
\begin{eqnarray*}\label{3.21}
&&t\log q(t,x,y)\nonumber\\
&\geq & t \log \inf_{z\in B(y,\varepsilon)}q(\eta t,z,y) + t \log \bP_x \left( X_{(1-\eta)t}\in B(y,\varepsilon) \right) \nonumber\\
&\geq & t\log \left( C_1 e^{-C_2t} (\eta t)^{-{d}/{2}}\right)
-\frac{C_3 \varepsilon^2}{2\eta }  + t\log \left( (2\pi (1-\eta)t)^{-{d}/{2}} \omega_d (\varepsilon-\delta)^d\right) \nonumber\\
& &  -\frac{(|x-y|+\varepsilon-\delta)^2}{2(1-\eta)}  +t\log \Big( 1-2  (2\pi (1-\eta)t)^{{d}/{2}}
\omega_d^{-1} (\varepsilon-\delta)^{-d} \nonumber\\
&& \ \times  \exp \Big( \frac{(M+\varepsilon-\delta)^2}{2(1-\eta)t}  -\frac{\delta^2}{8(1-\eta)t  C_4^2 e^{2 C_5}N^{C_6}_{(1-\eta)t} \big( \hbox{$\sum_{i=1}^d$} |\mu_i|\big)^2 }  \Big) \Big)
\end{eqnarray*}
Since $N^{C_6}_{(1-\eta)t} \left( \hbox{$\sum_{i=1}^d$} |\mu_i|\right)^2\rightarrow 0$ as $t\rightarrow 0$,
it follows that
\begin{eqnarray}\label{3.21}
\liminf_{t\rightarrow 0}t\log q(t,x,y)
&\geq & - \frac{C_3 \varepsilon^2}{2\eta}-\frac{(|x-y|+\varepsilon-\delta)^2}{2(1-\eta)}
\end{eqnarray}
uniformly in $x,y$ with $|x-y|\leq M$.
Now letting first $\varepsilon\rightarrow 0$ and then $\eta\rightarrow 0$,  we have
$$
\liminf_{t\rightarrow 0}t\log q(t,x,y)\geq -\frac{|x-y|^2}{2}
$$
uniform in $x, y\in \bR^d$ with $|x-y|\leq M$.  \qed

\medskip

Combining Propositions \ref{P:3.3} and \ref{P:3.4} establishes  Theorem \ref{T:3.1}.

\bigskip

\noindent{\bf Acknowledgement}. We thank James Norris and Renming Song for useful discussions. This work is partially supported by Simons Foundation Grant 520542,
  and NNSF of China (11671372,  11431014,  11401557,  11731009).

\end{document}